\documentclass[10.5pt]{article}

\usepackage{amsmath}
\usepackage{amssymb}
\usepackage{amsthm}
\usepackage{graphicx}
\usepackage{tikz}
\usepackage{xcolor}
\usepackage{dsfont}

\allowdisplaybreaks

\usepackage[margin=1.2in]{geometry}
\usepackage{bookmark}

\def\E{\mathbb{E}}
\def\var{\mathbb{Var}}

\def\E{\mathbb{E}}
\def\R{\mathbb{R}}

\def\Z{\mathbb{Z}}

\def\eps{\varepsilon}
\def\del{\delta}

\def\cF{\mathcal {F}}
\def\cG{\mathcal {G}}

\def\1{\mathbf{1}}

\def\F{\mathcal {F}}

\def\tce{t_c + \eps}
\def\tce2{t_c + \frac{\eps}{2}}
\def\ER{Erd\H{o}s-R\'{e}nyi }
\def\erdos{Erd\H{o}s }
\def\renyi{R\'{e}nyi }

\def\var{\text{var}}

\newtheorem*{theorem*}{Theorem}
\newtheorem{theorem}{Theorem}

\newtheorem{defn}[theorem]{Definition}
\newtheorem*{defn*}{Definition}
\newtheorem{prop}[theorem]{Proposition}
\newtheorem*{prop*}{Proposition}
\newtheorem{conj}[theorem]{Conjecture}
\newtheorem*{conj*}{Conjecture}

\newtheorem{question}{Question}
\newtheorem*{fact*}{Fact}

\newtheorem{task}{Task}

\begin{document}
\title{Searching for (sharp) thresholds in random structures: where are we now?}

\author{Will Perkins\footnote{School of Computer Science, Georgia Institute of Technology, Atlanta, GA, USA}}

\date{\today}
\maketitle

\begin{abstract}
We survey  the current state of affairs in the study of thresholds and sharp thresholds in random structures on the occasion of the recent proof of the Kahn--Kalai Conjecture by  Park and Pham and the fairly recent proof of the satisfiability conjecture for large $k$ by Ding, Sly, and Sun.  Random discrete structures appear as fundamental objects of study in many scientific and mathematical fields including statistical physics, combinatorics, algorithms and complexity, social choice theory, coding theory, and statistics.  While the models and properties of interest in these fields vary widely, much progress has been made through the development of general tools  applicable to large families of models and properties all at once.   Historically these tools originated   to solve or make progress on specific, difficult conjectures in the areas mentioned above.  We will survey recent progress on some of these hard problems  and describe some challenges for the future.

This survey was prepared in conjunction with a talk for the Current Events Bulletin at the 2024 Joint Mathematics Meetings in San Francisco.
\end{abstract}

\section{Introduction}
\label{secIntro}

Randomness is a powerful tool for algorithm design, scientific discovery, modeling the world,  engineering, and mathematical proof.   To use randomness effectively and confidently, we would like to accurately understand the properties of a typical random object or structure.  Perhaps most familiar is the use of the Central Limit Theorem to understand statistical significance and margin of error in scientific studies and political polling. 

In modern applications in algorithms, physics, social science, and other fields the random objects of interest can be very large and complex (think large computer networks, social networks, neural networks, or models of many interacting particles).  We would like to be able to understand with as much accuracy as possible, what properties random structures typically possess, and how these typical properties change as underlying parameters change. 

Even for the simplest class of random structures (those in which elements of a structured set are included independently with the same probability) this question can be extremely challenging, depending on the complexity of the property of interest.  Over the past 40 or more years, specific problems have driven different fields (including probabilistic combinatorics, algorithms, and statistical physics) to develop powerful tools for investigating this type of problem.  Some of these tools are designed for  specific settings while others are very general. 

Recently major progress on both general and specific problems has been made: Park and Pham proved the `Kahn--Kalai Conjecture'~\cite{park2023proof} and Ding, Sly, and Sun proved the `Satisfiability Conjecture' for large $k$~\cite{ding2022proof}.
The aim of this survey will be to put these exciting developments in a shared context and give an idea of some of the remaining and pressing challenges in the general area of thresholds in random structures.
A key theme will be  the way in which probability, combinatorics, and algorithms interact in these questions.

\subsection{Random structures, monotone properties, and thresholds}
\label{secProperties}

The main setting for this survey  will be $p$-biased product probability measures on the discrete cube $\{0,1\}^N$; that is, associating to a vector $x \in \{0,1\}^N$ the set $S_x= \{ i \in [N] : x_i=1\}$, probability measures of the form 
\[ \mu_p(S) = p^{|S|}(1-p)^{N - |S|} \,. \] 
Typically we will think of $N$ being large, and $p$ small, possibly depending on $N$.   This is simply a random subset with elements chosen independently with the same probability.  Often the set of coordinates $[ N ] = \{1, \dots, N\}$ will have additional structure; the coordinates might represent vertices or edges of some graph, or elements of an ordered set, or many other possibilities.  Then we can think of the random set $S \sim \mu_p$ as being a \textit{random structure}.   

We will see numerous examples shortly, but perhaps the first example to have  in mind  is flipping a sequence of $N$ independent coins each with probability $p$ of landing heads.  The structure here is the ordering of coordinates first to last.  

A \textit{property} $\cF$ is simply an event in this probability space, $\cF \subseteq \{0,1\}^N$.  A property $\cF$ is non-trivial if $\cF \ne \emptyset$, $\cF \ne \{0,1\}^N$.  In our example, non-trivial properties  include `Flipping all heads'; `Flipping at least $N/3$ heads'; `Flipping an odd number of heads', `Flipping a tail on the 7th and 11th flips', and so on.  A property is \textit{monotone} if it is closed under changing $0$'s to $1$'s; that is, if $x \in \cF$ and $y \ge x$ coordinate-wise, then $y$ must also be in $\cF$.   In the examples above about coin flips, all properties are monotone (or have a monotone complement) except the property of flipping an odd number of heads.

For a non-trivial monotone property, the probability that $\cF$ holds under $\mu_p$ increases from $0$ to $1$ as $p$ increases from $0$ to $1$, and we will be interested in how rapidly it does so.  Take the property of flipping at least $N/3$ heads.  By using a normal approximation to a binomial, we see that when $p = \frac{1}{3} + \frac{c \sqrt{2}}{3\sqrt{N}}$, $\mu_p(\cF)$ is approximately equal to the probability that a standard normal is at least $-c$; and so as $p$ passes $1/3$ in an interval of length on the order $N^{-1/2}$,  $\mu_p(\cF)$ jumps rapidly from near $0$ to near $1$.   The rest of this survey is devoted to understanding similar phenomena for much more complex properties of random structures.

\subsection{Examples}

Here we give some examples from combinatorics, probability, statistical physics, computer science, and social choice theory.  These  examples have driven considerable interest in thresholds in random structures from different fields and their study from these different perspectives has also given the field shared language and intuition.

\subsubsection{Probabilistic combinatorics}

\begin{figure}
\begin{minipage}[c]{0.44\linewidth}
\includegraphics[width=\linewidth]{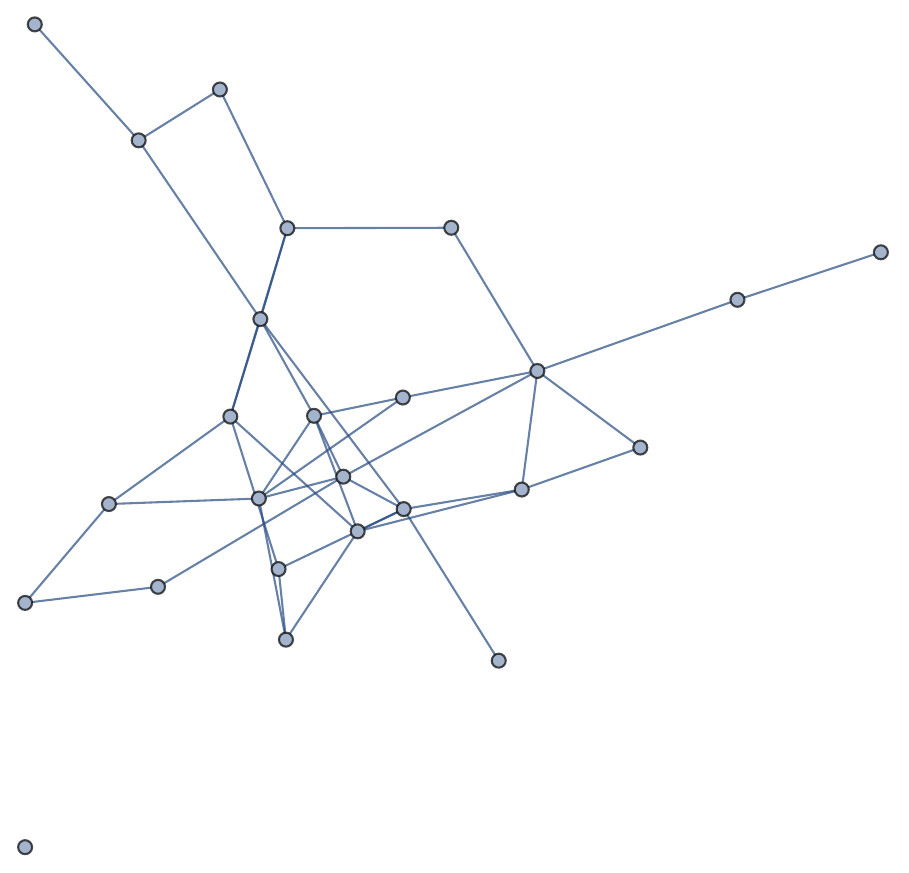}
\caption{A realization of the random graph  $G(25,1/8)$}
\label{figGnpL}
\end{minipage}
\hfill
\begin{minipage}[c]{0.44\linewidth}
\includegraphics[width=\linewidth]{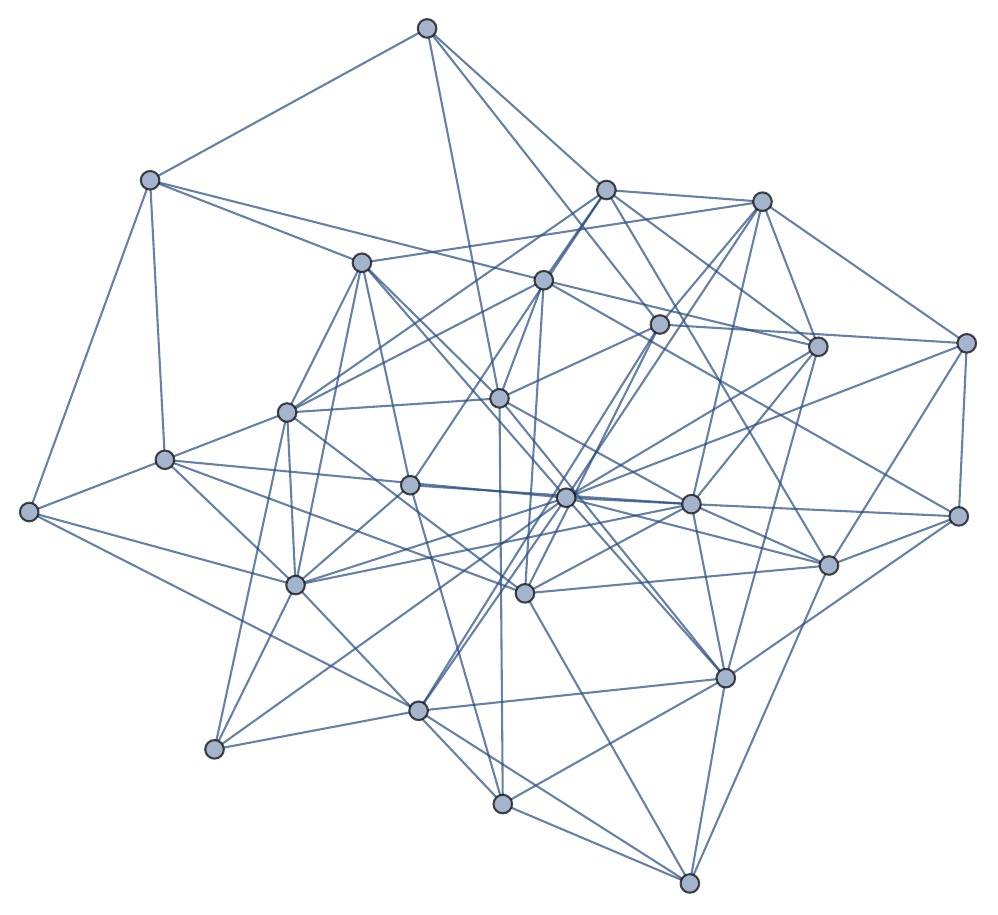}
\caption{A realization of the random graph  $G(25,1/4)$}
\label{figGnpH}
\end{minipage}\
\end{figure}

A central example comes from combinatorics.  Here we take $n \in  \mathbb N$ and let $N= \binom{n}{2}$, with coordinates representing edges of the complete graph $K_n$ on $n$ vertices.  Then the elements of $ \{0,1\}^N$ are in correspondence with (labeled) graphs on $n$ vertices.  A sample from $\mu_p$ is a random graph on $n$ vertices in which each possible edge is included independently with probability $p$; the model is known as the \textit{\ER random graph}~\cite{gilbert1959random,erdHos1960evolution,janson2011random,frieze2016introduction}.  See two samples from $G(n,p)$ in Figures~\ref{figGnpL} and~\ref{figGnpH}.   Despite the simplicity of its definition, the \ER random graph  exhibits a wide range of fascinating behaviors, and has been studied from many angles: as an interesting random structure in its own right, as a source of examples and counterexamples in graph theory, as a source of conjectured hard instances of combinatorial optimization problems, and as a host graph for random processes in probability and statistical physics.

\subsubsection{Computer science}

In computer science, a computational decision problem can be represented by a property of the discrete cube; or equivalently, \textit{boolean function} $f: \{0,1\}^N \to \{0,1\}$, where the input $x$ is encoded as a binary string and $f(x)=1$ if $x$ is a YES instance and $0$ otherwise.   A graph-theoretic decision problem corresponds directly to a property of graphs, as above.  The focus in computer science is on the computational resources (time, space, etc.) required to compute a given boolean function $f$, either in the worst case over inputs, or, as we discuss more below, for typical or random instances $x$, such as those drawn according to $\mu_p$.

\subsubsection{Statistical physics}

Much of the language and intuition around thresholds and sharp thresholds comes from physical systems that exhibit sudden changes in qualitative behavior as some parameter is changed in a small way.  This is the phenomenon of a \textit{phase transition}.  The phase transitions of water from liquid to solid as temperature drops below $32 ^ \circ  \text{ F}$ or from liquid to gas as temperature rises past $212 ^ \circ \text{ F}$  are familiar to everyone.

Mathematically, phase transitions are non-analytic points (or discontinuities in functions or their derivatives) of observables of infinite systems as some parameter varies.  This is not quite the setting described above of sequences of finite random structures with parameters that vanish as the system size grows, but there is much to be gained by pushing this analogy in both directions; in Section~\ref{secCSP} we describe the impact of the phase transition perspective on computer science, while~\cite{duminil2019sharp} takes the analogy in the other direction.

\begin{figure}
\begin{centering}
\includegraphics[width=.38\linewidth]{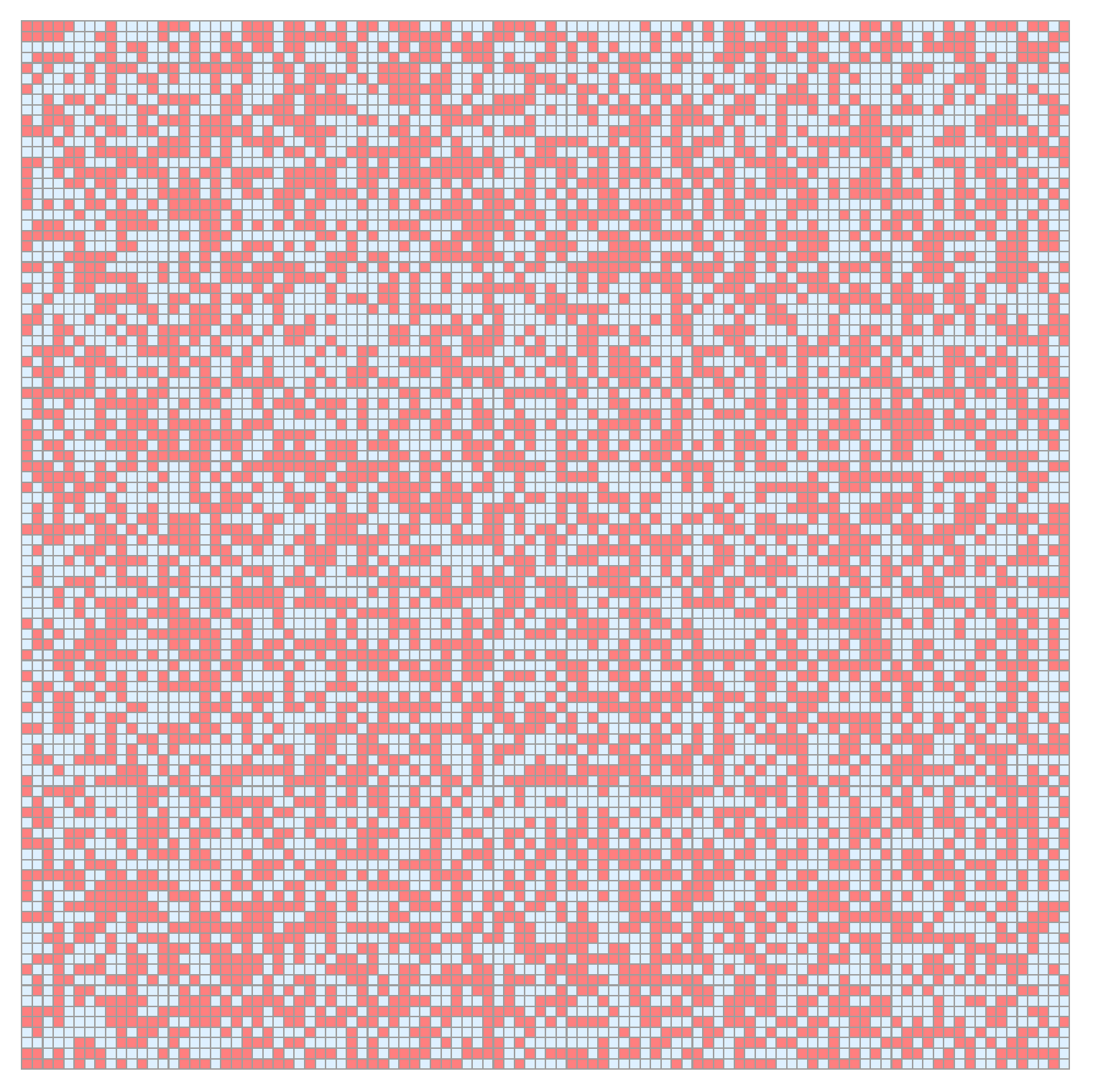}
\caption{Site percolation on a $100 \times 100$ grid with $p=1/2$}
\label{figPerc}
\end{centering}
\end{figure}

The statistical physics model perhaps most related to $\mu_p$ is that of \textit{percolation}~\cite{kesten1982percolation,grimmett1999percolation,bollobas2006percolation}.  Here the setting is typically an infinite graph, say a lattice like $\Z^d$ or the hexagonal lattice.  In site (resp. bond) percolation each vertex (resp. edge) is declared `open' independently with probability $p$; the main question is about the existence or non-existence of an infinite, connected `open' component.  See Figure~\ref{figPerc} for a depiction of site percolation on a finite portion of $\Z^2$.

\subsubsection{Social choice theory}

Social choice theory is the study of collective decision making, including the study of voting mechanisms.  We can view an election between two candidates with $N$ voters as a boolean function $f: \{0,1\}^N \to \{0,1\}$, where the candidates are labeled $0$ and $1$, and a vector $x \in \{0,1\}^N$ is the list of the votes of the $N$ voters.  The function $f$ is a rule for determining  the winner of the election given the votes.   A monotone function $f$ can be thought of as a rule with the sensible property that if a voter changes their vote from $0$ to $1$ the outcome of the election cannot change from $1$ to $0$.  Social choice theory lends a lot of evocative terminology to the study of boolean functions and thresholds.  A `dictator' boolean function is one that depends only on one coordinate; a `junta' is one that depends only on a small number of coordinates.  The majority function (i.e.\ equal representation democracy) has some extremal properties with applications in approximation algorithms (`Majority is stablest'~\cite{khot2007optimal,mossel2005noise}), while the `tribes' function, a hierarchical majority function from~\cite{ben1985collective}, inspired  conjectures on general properties of boolean functions that have proved essential to the study of thresholds.

\subsection{Outline}

 In Section~\ref{secPhenomena} we define critical probabilities, thresholds, sharp and coarse thresholds, and scaling windows, then illustrate these notions with three examples from the study of random graphs.  In the next three sections, we present some problems and conjectures that have driven progress on thresholds in different fields. In Section~\ref{secLatin} we present the first of these main motivating questions, the question of thresholds for spanning structures in random graphs and hypergraphs. We discuss perfect matchings in hypergraphs, Latin squares, and the Kahn--Kalai Conjecture.  In Section~\ref{secCSP} we discuss the random $k$-SAT model, the satisfiability conjecture, and  statistical physics predictions for thresholds in random computational problems.  In Section~\ref{secStat} we describe the stochastic block model and the predicted information-theoretic and algorithmic thresholds in the model.

We conclude in Section~\ref{secConclude} with some questions and open problems.

\subsection{Further reading}

This survey is not meant to be exhaustive at all, but rather aims to describe some recent results and open problems in different areas, all around the topic of thresholds and sharp thresholds in random structures.  This area is very fortunate to have several excellent books and surveys devoted to the topic from different perspectives.   Some of these sources are listed below with a few remarks.  

Friedgut's 2005 survey `Hunting for sharp thresholds'~\cite{friedgut2005hunting}, following the developments in~\cite{friedgut1999sharp,achlioptas1999sharp,friedgut2000sharp}, is devoted to the question of which properties have sharp thresholds (defined below), criteria for proving sharpness of a threshold, and some intuition and meta conjectures on the topic.  

Kalai and Safra's 2006 survey~\cite{kalai2006perspectives} takes a broad look at threshold phenomena in computer science, mathematics, and social choice theory and explains how notions arising in the latter (of `influence' or `pivotality' of coordinates) along with Fourier analysis of boolean functions can be used to study the sharpness or coarseness of thresholds.  O'Donnell's textbook on discrete Fourier analysis of boolean functions~\cite{o2014analysis} is a great reference to learn about these tools.

From a different perspective, Duminil-Copin's 2019 survey~\cite{duminil2019sharp} describes how general tools from the study of sharp thresholds of boolean functions can be used to study phase transitions in classical statistical physics models like percolation on infinite graphs.

Rao's recent survey~\cite{rao2023sunflowers} describes the 2021 breakthrough of Alweiss, Lovett, Wu, and Zhang on sunflowers~\cite{alweiss2021improved}, the method of which led to the proof of the fractional Kahn--Kalai conjecture by Frankston, Kahn, Narayanan, and Park~\cite{frankston2021thresholds}.  Even more recently Park has written an expository article~\cite{park2023threshold}  explaining the intuition behind the Kahn--Kalai conjecture and the consequences of its proof by Park and Pham~\cite{park2023proof}. 

Finally, in Sections~\ref{secCSP} and~\ref{secStat} on random computational problems and statistical inference we discuss areas that have been shaped to a great extent by the field of statistical physics through questions, conjectures, and methods.  The textbook of Mezard and Montanari~\cite{mezard2009information} is a great resource for understanding these methods, while the surveys of Zdeborov{\'a} and Krzakala~\cite{zdeborova2016statistical}, Moore~\cite{moore2017computer}, and Abbe~\cite{abbe2017community} give an account of  developments in this area in the last decade.

\section{General questions, notions, and phenomena}
\label{secPhenomena}

In what follows, the perspective will be asymptotic as some underlying parameter, usually $N$ or $n$, tends to infinity.  We use  standard asymptotic notation: $O( \cdot)$, $o(\cdot)$, $\Omega(\cdot)$, $\omega(\cdot)$, $\Theta(\cdot)$ to compare growth rates of functions.  We use a subscript, e.g. $O_\eps(\cdot)$, to indicate the implied constant may depend on $\eps$.  When discussing probabilities we say an event $A$ holds `with high probability' or `whp' as $N \to \infty$ if $\lim_{N \to \infty} \Pr(A) =1$; that is, $\Pr(A) = 1-o(1)$.  Here we implicitly  consider a sequence of probability spaces.

For a non-trivial monotone property $\cF$, the probability $\mu_p(\cF)$, as a function of $p$, is a strictly increasing function that starts with $\mu_0(\cF) = 0$ and $\mu_1(\cF)=1$ (in fact it is a polynomial in $p$).  The basic question about a monotone property of a random structure is how $\mu_p(\cF)$ increases from $0$ to $1$ as $p$ increases.  More generally we will think of a sequence of random structures indexed by $N$ along with a sequence of monotone properties, and ask how $\mu_p(\cF_N)$ increases from $0$ to $1$ as a function of $N$.   From here on, we will write $\cF$ and $\mu_p$ even though both implicitly depend on $N$.

Because $\mu_p(\cF)$ is strictly increasing, we can define the \textit{critical probability} of $\cF$, $p_c(\cF)$, as
\[ p_c(\cF) = \{ p : \mu_p(\cF) = 1/2 \} \,. \]

The first task in studying a particular property $\cF$ is to identify, at least approximately, the critical probability. 
\begin{task}
\label{taskPc}
For a given  monotone property $\cF$, determine $p_c(\cF)$ asymptotically as $N \to \infty$, or at least determine the asymptotic order of $p_c(\cF)$. 
\end{task}

Obtaining the asymptotics of $p_c$ means finding some $f(n)$ so that $p_c(\cF) = (1+o(1)) f(n)$; finding the asymptotic order means finding $f(n)$  so that  $p_c(\cF) = \Theta(f(n))$. 

In the context of random graphs, this task was put forward in the original paper of Erd\H{o}s and R\'{e}nyi~\cite{erdHos1960evolution} and has been a central topic in probabilistic combinatorics since.  In Section~\ref{secLatin} we discuss one of the most challenging classes of problems for this task: that of determining the asymptotic order of $p_c$ for the existence of different kinds of spanning subgraphs (including perfect matchings for random graphs and hypergraphs, $H$-factors of random graphs, and other combinatorial designs (including Latin squares) in random structures).  We will discuss the recent breakthrough work on the Kahn--Kalai conjecture in~\cite{frankston2021thresholds,park2023proof} that has led to a new powerful tool for determining the asymptotic order of $p_c$.

After  Task~\ref{taskPc}, the most fundamental question about monotone properties is how quickly $\mu_p(\cF)$ increases from near $0$ to near $1$.    To quantify this, we define the \textit{width of the scaling window} of a monotone property.  Let
\[ T_\eps(\cF) = \{ p_{1-\eps} - p_\eps : \mu_{p_{1-\eps}}(\cF) =1-\eps \text{ and } \mu_{p_\eps}(\cF)= \eps \} \,. \]
In words $T_{\eps}$ is the amount $p$ has to increase for the probability of $\cF$ to increase from $\eps$ to $1-\eps$.

Erd\H{o}s and R\'{e}nyi~\cite{erdHos1960evolution} defined a \textit{threshold function} for a monotone property as follows: $p_{\text{t}}(n)$ is a threshold function for a monotone property $\cF$ if the following hold:
\begin{enumerate}
\item If $p =\omega(p_{\text{t}}(n))$, then $\mu_p(\cF) = 1-o(1)$. 
\item If $p =o(p_{\text{t}}(n))$, then $\mu_p(\cF) =o(1)$. 
\end{enumerate}
Equivalently, $\cF$ has a threshold function (and we may take the function to be $p_c(\cF)$) if for every small $\eps>0$,  $T_{\eps}(\cF) = O_{\eps}(p_c(\cF))$. 

The following result of Bollob{\'a}s and Thomason justifies the abstract study of thresholds in monotone properties and the specific definition of $p_c(\cF)$ above.

\begin{theorem}[Bollob{\'a}s--Thomason~\cite{bollobas1987threshold}]
\label{thmBT}
Every non-trivial monotone property has a threshold function; moreover one can take $p_c(\cF)$ to be this threshold function.
\end{theorem}
The original proof of this theorem uses the Kruskal--Katona theorem, but a simple modern proof is as follows.  Let $K$ be a large positive integer.  If $\mu_p(\cF) =\eps$, then $\mu_{K p} (\cF) \ge 1- (1-\eps)^K$ by superimposing $K$ independent copies of the random structure and applying monotonicity; taking $K \approx \log (1/\eps)/\eps$ proves the theorem.

Erd\H{o}s and R\'{e}nyi further classified thresholds as \textit{sharp} or \textit{coarse} according to how $T_\eps(\cF)$ compares to $p_c(\cF)$.  We say $\cF$ has a \textit{sharp threshold} if for every $\eps>0$,  $T_\eps(\cF) = o_{\eps}(p_c(\cF))$; otherwise $\cF$ has a coarse threshold.   Equivalently, $\cF$ has a sharp threshold if for every $\eps>0$, 
\begin{enumerate}
\item If $p \ge (1+\eps) p_c(\cF)$, then $\mu_p(\cF) = 1-o(1)$.
\item If $p \le (1-\eps) p_c(\cF)$, then $\mu_p(\cF) =o(1)$.
\end{enumerate}

\begin{task}
\label{taskSharp}
Determine if a monotone property $\cF$ has a sharp or coarse threshold.  
\end{task}

In fact, in very general settings Task~\ref{taskSharp} has been solved, with the machinery of discrete Fourier analysis.    Roughly, as Friedgut describes in~\cite{friedgut2005hunting} and Kalai and Safra write in~\cite{kalai2006perspectives}, we should expect a sharp threshold unless there is reason to expect a coarse threshold; and the reason to expect a coarse threshold is if $\cF $ is essentially determined by a small number of coordinates or by the presence or absence of a small substructure.  This is made rigorous in increasing generality in Friedgut's theorem on sharp thresholds~\cite{friedgut1999sharp}, Bourgain's theorem in the appendix of~\cite{friedgut1999sharp}, and Hatami's result in~\cite{hatami2012structure}. We describe in Section~\ref{secCSP} the random $k$-SAT model, a main motivation for Friedgut's important result.

A final task in understanding a particular property $\cF$ with a sharp threshold is getting more precise bounds on the scaling window. 

\begin{task}
\label{taskScaling}
Determine the asymptotic order of the width of the scaling window of a monotone property $\cF$. 
\end{task}

Unlike with Tasks~\ref{taskPc} and~\ref{taskSharp}, as of yet there is no general principle in determining the magnitude of $T_{\eps}(\cF)$ beyond distinguishing sharp from coarse: the width of the scaling window really seems to depend on the particular details of the property $\cF$ in question, and being able to determine the asymptotic order of the width indicates a near-complete understanding of the property.

\subsection{Examples}
\label{sec3examples}

We give three  examples to illustrate these different tasks. The examples are the properties in the random graph $G(n,p)$ of \textit{containing a triangle}, \textit{being connected}, and \textit{being 3-colorable}; call these properties $\cF_{K_3}$, $\cF_{\text{connected}}$, $\cF_{3-\text{col}}$ respectively. 

\begin{figure}
\begin{minipage}[c]{0.48\linewidth}
\includegraphics[width=\linewidth]{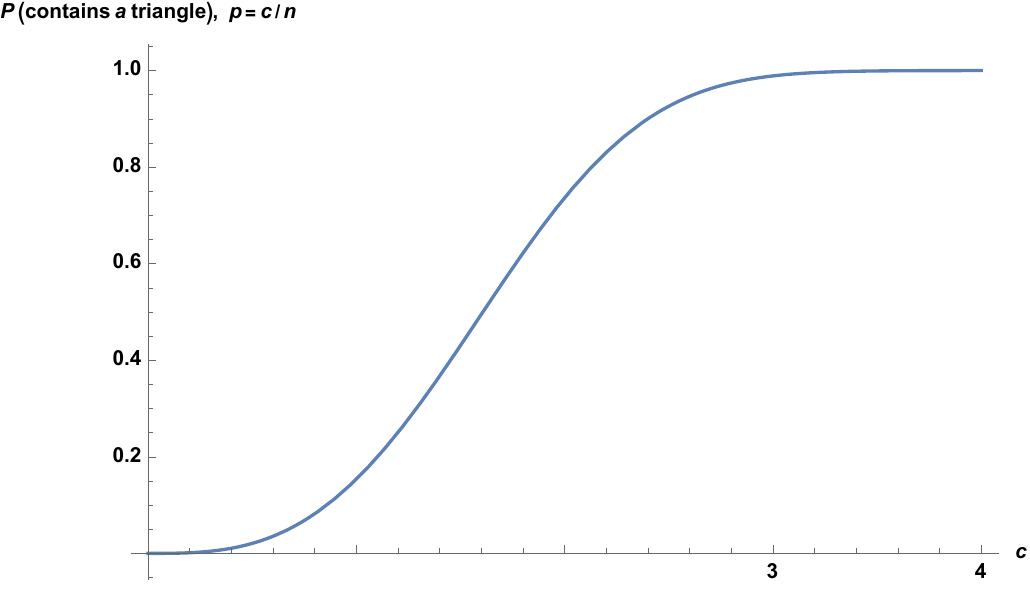}
\caption{Scaling window of $\cF_{K_3}$}
\label{figTri}
\end{minipage}
\hfill
\begin{minipage}[c]{0.48\linewidth}
\includegraphics[width=\linewidth]{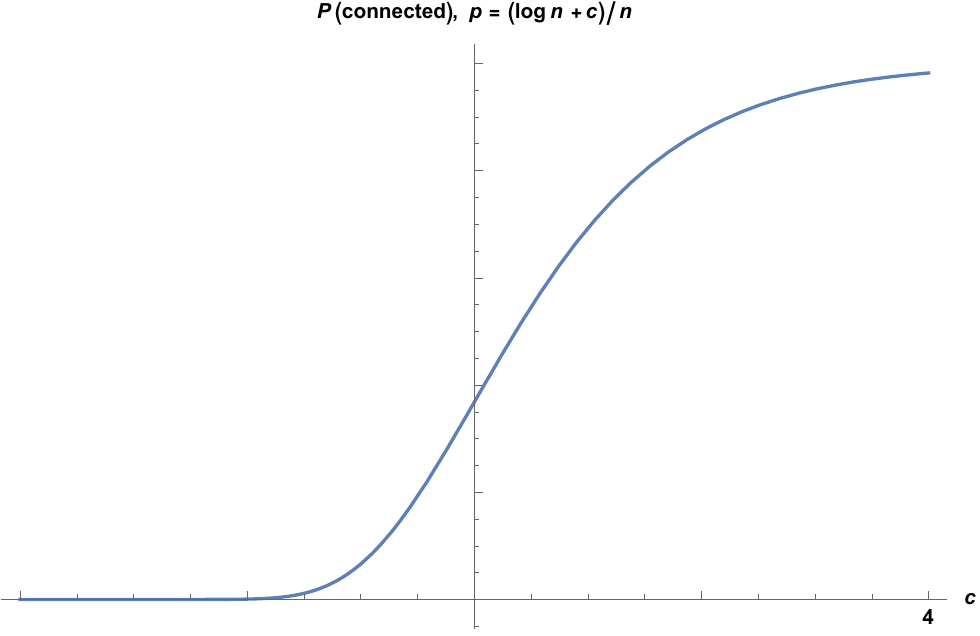}
\caption{Scaling window of $\cF_{\text{connected}}$}
\label{figCon}
\end{minipage}
\end{figure}

To understand the probability that $G(n,p)$ contains a triangle, we can use the first- and second-moment methods (see e.g.~\cite{alon2016probabilistic} for an exposition).  Let $X$ be the number of triangles.  Then $\E X = \binom{n}{3} p^3$ and $\var(X) = (1+o(1))\binom{n}{3} p^3   $.   Using this we can determine that $p_c(\cF_{K_3}) = \Theta(1/n)$.  Markov's Inequality yields:
\begin{align*}
\Pr(\cF_{K_3}) &= \Pr(X \ge1 ) \le \E X = (1+o(1)) \frac{n^3 p^3}{6} = o(1)
\end{align*}
when $p = o(1/n)$. Chebyshev's Inequality yields
\begin{align*}
\Pr( \overline{\cF_{K_3}})&= \Pr(X=0) \le \Pr (|X- \E X| \ge \E X) \le \frac{\var(X)} {(\E X)^2}= (1+o(1)) \frac{6}{n^3p^3} =o(1)
\end{align*}
when $p=\omega(1/n)$.  This pair of facts tells us that $p= c/n$ (for any constant $c>0$) is a threshold function for 
$\cF_{K_3}$.   In fact much more is known:  when $p= c/n$, the distribution of $X$ converges to a $\text{Poisson}(c^3/6)$ random variable. Thus  the scaling window is of length $\Theta(1/n)$ and  the threshold is coarse; see Figure~\ref{figTri}.  

Next we turn to the property $\cF_{\text{connected}}$.  Being connected is a more `global' property  than that of containing a triangle, but there is an easy  lower bound  on $p_c(\cF_{\text{connected}})$: if $G$ has an isolated vertex (and $n >1$) then $G$ must be disconnected. Or, taking complements, $\cF_{\text{connected}} \subseteq \cF_{\text{no isolated vertex}}$. As with triangles it is straightforward to understand the threshold for containing isolated vertices  using the method of moments.  Let $Y$ be the number of isolated vertices; then   $\E Y = n (1-p)^{n-1}$.  As with triangles, when $\E Y$ tends to a positive constant, the distribution of $Y$ converges to a Poisson random variable; when $\E Y \to 0$ whp $Y =0$; and when $\E Y \to \infty$, whp $Y \ge 1$. This tells us the exact scaling window for $\cF_{\text{no isolated vertex}}$: if $p = \frac{ \log n  + c }{ n}$ with $c \in \R$ constant, then $\Pr(\cF_{\text{no isolated vertex}}) = (1+o(1)) e^{- e^{-c}}$, see Figure~\ref{figCon}.  In this case the width of the scaling window is $\Theta(1/n) = o(p_c)$ and so the threshold is sharp.

It turns out that a first-moment argument on connected components of size $\ge 2$ from  \erdos and \renyi\cite{erdds1959random} shows that in fact $\cF_{\text{no isolated vertex}}$ approximates $\cF_{\text{connected}} $ very well indeed around the threshold.  An effective way to state the approximation is as a \textit{hitting time result}, which we now describe.

We can couple the random graphs $G(n,p)$ for all $p \in [0,1]$ simultaneously by drawing iid Uniform$[0,1]$ random variables $U_{ij}$, $1 \le i < j \le n$. Then we form $G(n,p)$ by including $(ij) \in E$ if and only if $U_{ij} \le p$. Under this coupling $G(n,p')$ is a subgraph of $G(n,p)$ if $p' \le p$.   We can then define the `random graph process' as the discrete time process $G_0, G_1, \dots, G_{\binom{n}{2}}$ where at each step a uniformly random edge (not already present) is added to the graph.  This process is what results from raising $p$ from $0$ to $1$ and marking the appearance of new edges in the coupling. 

\begin{theorem}[Hitting time for connectivity~\cite{erdds1959random}]
Whp over the random graph process we have
\[ \min_t \{ G_t \text{ is  connected}  \}= \min_t \{ G_t \text{ has no isolated vertices}\} \,.  \]
\end{theorem}
That is, the  edge that touches the last isolated vertex in the process also connects the graph, with high probability over the process. This implies  that the threshold and scaling widow are the same for $\cF_{\text{connected}}$ and $\cF_{\text{no isolated vertex}}$. 

We will see below in Section~\ref{secLatin} further examples of hitting time results along the same lines for  more complicated properties than connectivity.

Finally, consider a property for which we do not yet have a full understanding.  Let $\cF_{3-\text{col}}$ be the property that $G$ has a proper $3$-coloring of its vertices (no monochromatic edges).  The complement of this property is a non-trivial monotone property.  

One can show that if  $p=c/n$, $c<1$, then whp $G(n,p)$ is $3$-colorable by showing that whp all connected components  are trees or unicyclic (and thus $3$-colorable).  For an upper bound on $p_c$ one can use the first moment method on $Z$, the number of $3$-colorings of $G$.  As a  first attempt we can bound $\E Z \le 3^n (1-p)^{\frac{n^2}{6}}$ using the fact that a balanced partition of $n$ vertices has the fewest potential monochromatic edges.  This tells us that $G(n,c/n)$ is not $3$-colorable whp when  $c> 6 \log 3 \approx 6.5917$.   However, we can do something a little  bit more clever: whp $G(n,c/n)$ has $\frac{cn}{2} +o(n)$ edges and so we can condition on the number of edges $m$, and compute $\E Z$ in the random graph $G_m$ from the random graph process.  This gives the bound 
\begin{align*}
\E_{G_m} Z \le 3^n \left( 1 - \frac{1}{3} \right)^m 
\end{align*}
which tends to $0$ when $m \sim cn/2$, $c > \frac{2 \log3 }{ \log 3 - \log 2} \approx 5.419$.  We will see more of this idea of combining the moment method with conditioning on typical events in Section~\ref{secCSP}.

Together these facts tell us that $p_c(\cF_{3-\text{col}}) = \Theta(1/n)$ and the width of the scaling window is $O(1/n)$, uniformly over $\eps$.  What are the asymptotics of $p_c$?  Is there some $d_3$ so that $G(n,d/n)$ is $3$-colorable whp when $d< d_3$ and whp not $3$-colorable when $d> d_3$?    In fact this is the only remaining open problem from the original \erdos and \renyi paper on random graphs~\cite{erdHos1960evolution}.  In Section~\ref{secCSP} we describe   predictions for this problem  made by statistical physicists.

\section{Perfect matchings, Latin squares, Kahn--Kalai, and spread distributions}
\label{secLatin}

From the beginning of the study of random graphs, a key question has been `For a given subgraph $H$, when in the evolution of the random graph does a copy of $H$ first appear?' That is, what is the threshold for the property $\cF_H$  of containing  a copy of $H$.

For fixed size $H$, the answer is now known completely~\cite{bollobas1998random}, though more delicate in general than the case of $H= K_3$ above. A threshold  function for  $\cF_H$ is   the \textit{expectation threshold} for the densest subgraph $H'$ of $H$ (maximizing the ratio of edges to vertices over all induced subgraphs of $H$); the value  $p_E(H')$ so that the expected number of copies of $H'$ in $G(n,p_E(H'))$ is $1$.  For $p = o( p_E(H'))$, the expected number of copies of $H'$ tends to $0$, and so whp $G(n,p)$ has no copies of $H'$ and thus no copies of $H$. When the expected number of copies of $H'$ for each $H' \subseteq H$ tends to infinity,  it has been showed that whp a copy of $H$ exists.

It can be much more difficult to determine the threshold for the appearance of subgraphs whose size grows with $n$; in particular, for \textit{spanning subgraphs} that include each of the $n$ vertices.  Two classical examples are those of $H$ being a perfect matching (with $n$ even) and $H$ being a Hamilton cycle.  As with connectivity, there are natural lower bounds for $p_c(\cF_H)$ in these cases, coming from local obstructions.  If $G$ has an isolated vertex, it cannot have a perfect matching; if $G$ has a vertex of degree  at most $1$ it cannot have a Hamilton cycle.  
In both cases these obstructions determine the thresholds, even in the strong form of hitting-time results~\cite{erdHos1966existence,komlos1983limit,bollobas1985random} that give a complete understanding of the scaling windows.

Other spanning subgraph problems remained open for many years.  One famous such problem is Shamir's Problem~\cite{schmidt1983threshold} (mentioned first in print by \erdos in~\cite{erdHos1981combinatorial}): the problem of determining the threshold for a random $k$-uniform hypergraph to contain a perfect matching.

A $k$-uniform hypergraph $G=(V,E)$ is a set of vertices $V$ along with a collection $E$ of $k$-sets of $V$ (the case $k=2$ is a usual graph).  The random $k$-uniform hypergraph $G^{(k)}(n,p)$ is a $k$-uniform hypergraph on $n$ vertices in which each possible $k$-set is an edge independently with probability $p$.   A perfect matching in a $k$-uniform hypergraph on $n$ vertices is a set of $n/k$ vertex-disjoint hyperedges that (necessarily) cover all $n$ vertices; $n$ must be divisible by $k$ for a perfect matching to be possible.  For $k=2$ (a graph) this coincides with the usual definition of a perfect matching.    

As with perfect matchings in graphs, the natural obstruction to a perfect matching in a hypergraph is a vertex not appearing in any hyperedges; this immediately yields a lower bound of $\Omega(\log n \cdot n^{1-k})$ for the threshold for perfect matchings.  

After many years of effort, Johansson, Kahn, and Vu proved that the lower bound is tight up to constants.

\begin{theorem}[Johansson--Kahn--Vu~\cite{johansson2008factors}]
\label{thmJKV}
Let $\cF$ be the property of $G^{(k)}(n,p)$ containing a perfect matching when $k$ divides $n$.  Then
\[ p_c(\cF) = \Theta_k( \log n \cdot n^{1-k}) \,. \]
\end{theorem}
Much more recently Kahn proved much finer results: first-order asymptotics of $p_c(\cF)$ in~\cite{kahn2023asymptotics} and the definitive hitting-time result~\cite{kahn2022hitting} (see also the related~\cite{heckel2021random,riordan2022random,heckel2023hitting}).

\subsection{The Kahn--Kalai conjectures}

Motivated by challenges like Shamir's Problem, Kahn and Kalai~\cite{kahn2007thresholds} made two  bold conjectures about thresholds.  These conjectures arise from the fact that it is often straightforward to give a lower bound for $p_c(\cF)$; that is, to show that when $p$ is small enough $\mu_p(\cF) =o(1)$. A matching or near matching upper bound is often much more difficult.  This is especially apparent in the examples above (connectivity, perfect matchings, etc.) in which there is a trivial local obstruction to some property.  Kahn and Kalai's two conjectures both have a similar flavor: `The best possible easy lower bound on $p_c$ is not too far from the truth.'  We present the two conjectures in the opposite order they appear in~\cite{kahn2007thresholds} (but in the order they arrived at them).

Their second conjecture is about the property $\cF_H$ of containing an isomorphic copy of some subgraph $H$ in a random graph (or hypergraph).  We have seen several examples of this for different $H$, including  a triangle,  a perfect matching, a Hamilton cycle, and  a hypergraph perfect matching.   As we saw for fixed-size subgraphs, it can be crucial to consider the rarest (in terms of expected number of copies) subgraphs of $H$.  This motivates the following definition of the \textit{subgraph expectation threshold} of $H$:
\begin{align*}
p_E(H) = \min \{ p : \E_{G(n,p)} X_{H'} \ge 1 \, \forall \text{ subgraphs } H' \subseteq H  \,\} \,.
\end{align*}
Clearly $p_E(H)$  bounds below the threshold of $\cF_H$; if $p = o(p_E(H))$ then there is some subgraph $H' \subseteq H$ that whp does not appear in $G(n,p)$. Kahn and Kalai conjecture this cannot be too far from the truth.
\begin{conj}[Kahn and Kalai~\cite{kahn2007thresholds}]
\label{conjKKsubgraph}
There exists an absolute constant $K$ so that for all graphs $H$, 
\[ p_c(\cF_H) \le  K \cdot p_E(H) \cdot \log  |V(H)| \,. \]
\end{conj}

For $H$ of fixed size, Conjecture~\ref{conjKKsubgraph} asserts that $p_e(H)$ is a threshold for $\cF_H$ and this is known to be true as described above.  The real content of the conjecture is for $H$ that grow with $n$.   In particular Conjecture~\ref{conjKKsubgraph} is tight (up to the constant $K$) for perfect matchings in graphs and  hypergraphs.

Their first conjecture was in the more general setting of monotone properties of $\{0,1\}^N$.  To state it we must define a more abstract notion of an expectation threshold (following~\cite{talagrand1995all,talagrand2006selector,kahn2007thresholds}).  A property $\cF \subseteq \{0,1\}^N$ is $p$-\textit{small} if there exists a `cover' $\cG \subseteq \{0,1\}^N$ so that:
\begin{enumerate}
\item $\forall \, T \in \cF \,  \,  \exists \, S \in \cG, S \subseteq T $.
\item $\sum_{ S \in \cG} p^{|S|} \le \frac{1}{2} $.
\end{enumerate}
To give some intuition for this definition, if $H$ is a graph and  $H' \subseteq H$, then $\F_{H'}$ is a cover of $\cF_H$, satisfying the first condition above.

The \textit{expectation threshold}  of $\cF$, $q (\cF)$, is the largest $p$ for which $\cF$ is $p$-small.  Finally let $L(\cF)$ be the maximum size of a minimal element of $\cF$.

Notice that $q(\cF)$ is a lower bound for $p_c(\cF)$. Let $q=q(\cF)$. Then
\begin{align*}
\mu_q(\cF) &\le \sum_{S \in \cG}  \, \,\sum _{T\in \cF, T \supseteq S} \mu_q(T)  \\
&\le  \sum_{S \in \cG}  \, \,\sum _{T \supseteq S} \mu_q(T)  \\
&= \sum_{S \in \cG} q^{|S|} \le \frac{1}{2} \,.
\end{align*}
Kahn and Kalai conjectured that in the abstract setting this cannot be far from the truth for $p_c$.  Park and Pham then proved this.
\begin{theorem}[Kahn--Kalai Conjecture, now Park--Pham Theorem~\cite{park2023proof}]
\label{thmKK}
There exists an absolute constant $K$ so that for every monotone $\cF$, 
\[ p_c(\cF) \le K \cdot q(\cF) \cdot \log L(\cF) \, .\]
\end{theorem}

Theorem~\ref{thmKK} gives a powerful method for bounding $p_c$ from above: find a good (or the best) $\cG$ with the properties above, and this will determine $p_c$ up to a $\log N$ factor or better. Looking ahead, we can characterize $q(\cF)$ via the following integer optimization problem over the variables $g(x)$, $x\in \{0,1\}^N$ which we interpret as the indicator vector of  a potential cover $\cG$ witnessing $\cF$ being $p$-small.
\begin{align*}
V(\cF,p) = &\min   \sum_{x \in \{0,1\}^N}  p^{|x|} g(x)  \\
&\text{subject to } \\
& g(x) \in \{0,1\} \text{ for all } x \in \{0,1\}^N \\
&\sum _{ S \subseteq T} g(x_S) \ge 1 \text{ for all } T \in \cF
\end{align*}
In particular, if $V(\cF,p) \ge 1/2$, then $q(\cF) \le p$, and by Theorem~\ref{thmKK}, $p_c(\cF) \le K \cdot p \cdot \log L(\cF)$. However, solving an integer program over $2^N$ variables (for all large $N$) seems somewhat formidable. 

\subsection{Duality and spread distributions}

Before Park and Pham proved Theorem~\ref{thmKK}, Frankston,  Kahn,  Narayanan, and Park~\cite{frankston2021thresholds} (following the breakthrough of Alweiss, Lovett, Wu, and Zhang on the sunflower conjecture~\cite{alweiss2021improved}) proved a related but weaker conjecture of Talagrand from~\cite{talagrand2010many}, dubbed the `Fractional Kahn--Kalai Conjecture'.  To state this, we first relax the integrality constraint above and define a linear program (in $2^N$ variables)
\begin{align*}
V_f(\cF,p) = &\min   \sum_{x \in \{0,1\}^N}  p^{|x|} g(x)  \\
&\text{subject to } \\
& g(x) \in [0,1] \text{ for all } x \in \{0,1\}^N \\
&\sum _{ S \subseteq T} g(x_S) \ge 1 \text{ for all } T \in \cF \,.
\end{align*}
Then define $q_f(\cF)$, the \textit{fractional expectation threshold} of $\cF$, to be the largest $p$ for which $V_f(\cF,p)\le 1/2$. Again it is easy to see that $p_c(\cF) \ge q_f(\cF)$, with the same proof as above.

Moreover, we have the relations
\begin{align*}
q(\cF) \le q_f(\cF) \le p_c(\cF) \le K \cdot q(\cF) \cdot \log L(\cF)  \le K \cdot q_f(\cF) \cdot \log L(\cF) \,, 
\end{align*}
where the second-to-last inequality is Theorem~\ref{thmKK}; the weaker inequality $p_c(\cF) \le K \cdot q_f(\cF) \cdot \log L(\cF)$ is the main result of~\cite{frankston2021thresholds}. 

\begin{theorem}[Frankston,  Kahn,  Narayanan, and Park~\cite{frankston2021thresholds}]
\label{thmFrac}
There exists an absolute constant $K$ so that for every monotone $\cF$, 
\[ p_c(\cF) \le K \cdot q_f(\cF) \cdot \log L(\cF) \, .\]
\end{theorem}
To bound $p_c$ with Theorem~\ref{thmFrac} it suffices to lower bound $V_f(\cF,p)$: if $V_f(\cF,p) \ge1/2$ then $p_c(\cF) \le K \cdot p \cdot \log L(\cF)$. 
The nice thing about the linear programing formulation is that we can use duality to give a lower bound.  We can write the dual as a linear program with variables $\nu(T)$ for each $T \in \cF$ and constraints for each $S \subseteq \{0,1\}^N$. 
\begin{align*}
V_f(\cF,p) = &\max   \sum_{T \in \cF}  \nu(T)  \\
&\text{subject to } \\
& \nu(T) \ge 0 \text{ for all } T \in \cF \\
&\sum _{  T \supseteq S} \nu(T) \le p^{|S|} \text{ for all } S \subseteq \{0,1\}^N \,.
\end{align*}
We are interested in showing $V_f(\cF,p) \ge 1/2$ to bound $p_c$ and via duality this can be accomplished by exhibiting a good $\nu(\cdot)$.  A slightly more elegant formulation is due to Talagrand~\cite{talagrand2010many}, who made the following definition of a \textit{spread probability distribution}.

\begin{defn}
Let $\cF \subseteq \{0,1\}^N$. A probability measure $\nu$ supported on $\cF$ is $p$-spread  if for all $S \subseteq \{0,1\}^N$, 
\[ \sum _{  T \supseteq S} \nu(T) \le2 p^{|S|} \,. \]
\end{defn}
Putting all of the above together we obtain a very useful theorem.
\begin{theorem}
\label{thmSpread}
There is an absolute constant $K$ so that the following is true. Let $\cF$ be a monotone property that supports a $p$-spread probability measure $\nu$.  Then
\[ p_c(\cF) \le K \cdot p \cdot \log L(\cF) \,.\]
\end{theorem}
Theorem~\ref{thmSpread} follows from Theorem~\ref{thmFrac}, duality, and the extra factor $2$ allowing us to take $\nu$ to be a probability measure.

\subsection{Applications}

Theorem~\ref{thmSpread} is a brand-new tool for proving upper bounds on $p_c(\cF)$ and we will recount here some spectacular applications. But first as a warm-up we will see how it can be used to establish $p_c(\cF_{\text{perfect matching}}) = \Theta\left(\frac{ \log n }{n} \right)$.  

Let $n$ be even, and let $\nu$  be the uniform distribution on perfect matchings of the complete graph $K_n$.   Let $A$ be some set of edges of $K_n$ and let $M \sim \nu$.  We want to bound the probability that $A \subseteq M$.  First note that if $A$ is not a matching then this probability is $0$.  Now suppose $A$ is a matching of size $k$. Then, letting $\text{pm}(G)$ denote the number of perfect matchings of $G$ and  $(a)_b = \frac{a!}{(a-b)!}$,
\begin{align*}
\Pr(A \subseteq M) = \frac{ \text{pm}(K_{n-2k})}{\text{pm}(K_n)} = \frac{(n-2k)! }   {2^{n/2-k} (n/2-k)!}   \frac {2^{n/2} (n/2)!}{n!} = \frac{2^k (n/2)_k} { (n)_{2k}}  \le \left( \frac{e }{n } \right)^k \,,
\end{align*}
and so $\nu$ is $\frac{e}{n}$-spread.  Applying Theorem~\ref{thmSpread} gives the result, though note that this is weaker than the sharp threshold obtained in~\cite{erdHos1966existence}.   Here spreadness is proved by counting, and it helps a lot that we have an explicit formula for $\text{pm}(K_n)$ to use.

From Theorem~\ref{thmSpread} and similar counting arguments (see e.g.~\cite{frankston2021thresholds}), one can also rather easily obtain  the asymptotic order of the threshold for perfect matchings in hypergraphs and $K_r$-factors as well as the threshold for containing any bounded-degree spanning tree, previously known through the long proofs of Johansson, Kahn, and Vu~\cite{johansson2008factors} and Montgomery~\cite{montgomery2019spanning} respectively.

  Perhaps even more exciting are applications which before the theorem were completely out of reach but which now can be approached through finding spread distributions.  We describe one application here but see also~\cite{pohoata2022sharp,kelly2023optimal,anastos2023robust} for other applications of either Theorem~\ref{thmSpread} or the notion of spread distributions.

A \textit{Latin square} of order $n$ is an $n \times n$ matrix with entries in $\{1, \dots, n\}$ in which no row or column contains repeated entries.  A Latin square of order $8$ is given in Table~\ref{figLatin}.
\begin{table}[htp]
\begin{center}
\begin{tabular}{|c|c|c|c|c|c|c|c|}
\hline
1& 3& 4& 8& 2& 6& 7& 5\\
\hline
6& 1& 2& 5 &7& 4& 8& 3\\
\hline
5& 2& 3& 1& 8& 7& 4& 6\\
\hline
4& 6& 5& 3& 1& 8& 2& 7\\
\hline
2& 4& 6& 7& 3& 5& 1& 8 \\
\hline
 7& 8& 1& 2& 6& 3& 5& 4 \\
 \hline
 8&5& 7& 6& 4& 1& 3& 2 \\
 \hline
 3& 7& 8& 4& 5& 2& 6& 1 \\
 \hline
\end{tabular}
\end{center}
\caption{A Latin square of order $8$}
\label{figLatin}
\end{table}%
One can ask about the existence of Latin squares of order $n$, the number of Latin squares of order $n$, or the existence of Latin squares with additional properties.  Johansson~\cite{johansson2006triangle} asked a very natural threshold question about Latin squares: if each position $(ij)$ in an $n\times n$ matrix is assigned a list $L_{ij}$ of allowed symbols in $\{1, \dots, n\}$ by including each element in $L_{ij}$ independently with probability $p$, what is the threshold in $p$ for the existence of a Latin square in which each entry appears in the corresponding list $L_{ij}$?  The lower bound comes from a simple local obstruction: for $p \le (1-\eps) \frac{\log n}{n}$ there will be empty lists whp.    A natural conjecture given the results and conjectures on perfect matchings and factors in random graphs is that this local obstruction determines the threshold: its asymptotic order, first-order asymptotics, and even in the form of a hitting time.  This conjecture is stated explicitly in~\cite{luria2019threshold}, and see also the closely related conjectures in~\cite{casselgren2016coloring,simkin2017left,kang2023thresholds,sah2023threshold}.
\begin{conj}
\label{conjLatin}
The property  $\cF_{\text{Latin}}$ has a sharp threshold at $p = \frac{\log n}{n}$
\end{conj}
 Before Theorem~\ref{thmSpread}, there was really no effective way to to prove upper bounds on $p_c(\cF_{\text{Latin}})$: the best known result was that $p_c\le1 -\del$ for some fixed $\del>0$~\cite{andren2013avoiding}.

Using Theorem~\ref{thmSpread} and finding a sufficiently spread distribution on Latin squares, Sah, Sawhney and Simkin proved Conjecture~\ref{conjLatin} up to subpolynomial factors: $p_c = n^{-1+o(1)}$~\cite{sah2023threshold}.  Intriguingly, this bound already established that the threshold is sharp without determining the threshold, by applying Friedgut's theorem discussed in the next section.  Following this, Kang, Kelly, K{\"u}hn,  Methuku, and Osthus~\cite{kang2023thresholds} proved $p_c = O(\log^2 n/n)$, within a factor $\log n$ of the conjecture.  Then very recently Jain and Pham, and independently Keevash,  established the correct order of the threshold.
\begin{theorem}[Jain--Pham~\cite{jain2022optimal}; Keevash~\cite{keevash2022optimal}]
\label{thmLatin}
\[ p_c( \cF_{\text{Latin}}) = \Theta \left( \frac{\log n}{n} \right) \,.\]
\end{theorem}
Using Theorem~\ref{thmSpread},  it suffices to construct an $O(1/n)$-spread distribution on Latin squares.  It is natural to expect that (as in the case of perfect matchings) the \textit{uniform distribution} is $O(1/n)$ spread, but unlike perfect matchings it is very challenging to enumerate Latin squares.  So instead, the above authors have constructed non-uniform  spread distributions using sophisticated tools from probabilistic and extremal combinatorics, namely iterative absorption (in \cite{sah2023threshold,kang2023thresholds}), analyzing the  Lov{\'a}sz Local Lemma probability distribution (in~\cite{jain2022optimal}), and analyzing a random greedy stochastic process (in~\cite{keevash2022optimal}).

\section{Random $k$-SAT and the satisfiability conjecture}
\label{secCSP}

The  most important open problem in theoretical computer science (and one of the most important in  all of mathematics) is the P vs NP question: can computational decision problems with polynomial-time checkable certificates be solved in polynomial time?  

Most computer scientists and mathematicians believe that P$\ne$NP and that many hard computational problems exist; these include classic problems like MAX-CUT, Max Independent Set, traveling salesman, boolean satisfiability, graph coloring, among many others.  It is believed that these problem are computationally intractable in the worst-case over instances.  On the other hand, large instances of these problems are solved every day (see for example a survey on the success of SAT solvers in~\cite{GaneshV20} or exact solutions to very large instances of the traveling salesman problem in e.g.~\cite{applegate1998solution,chvatal2010solution}).  What can explain this discrepancy?  One possibility is that `typical' instances of certain NP-hard problems are tractable while hard instances are exceptional.

This is one motivation for the study of \textit{average-case complexity}: the computational tractability or intractability of random instances of computational problems.
Average-case complexity is a huge topic with many fascinating facets and many mysteries (see e.g.~\cite{bogdanov2006average,wigderson2006p}).  Here we will deal with a small slice of the topic.

Recall that a boolean CNF formula is the AND of OR's of  literals (boolean variables and their negations).  A $k$-CNF formula is the AND of clauses of $k$ literals each.  For example, the following is a small $3$-SAT formula:
\[ ( x_7 \vee \overline x_2 \vee \overline x_3) \wedge  ( \overline x_1 \vee  x_8 \vee \overline x_4) \wedge ( x_2 \vee  x_4 \vee  x_1) \,.\]
The $k$-SAT decision problem is to determine if there exists an assignment of True and False to the boolean variables so that the given $k$-CNF formula evaluates to true; if such an assignment exists the formula is satisfiable and unsatisfiable otherwise.

An empirical observation about algorithms brought the study of random computational problems to the attention of computer scientists, statistical physicists, and probabilists.   In 1992, Mitchell, Selman, and Levesque~\cite{mitchell1992hard} (working on problems in Artificial Intelligence and computational deductive reasoning)  generated uniformly random 3-SAT instances on $n$ variables with $m$ constraints, for fairly large $n$ and different values of $m$.  They ran standard heuristic SAT solving and SAT refutation algorithms and observed the following: the running time required to find a solution (or find a proof that none existed) showed a sharp peak (as a function of $m$) right around the point at which an estimate of the probability of such a random instance being satisfiable made a sharp decrease from near $1$ to near $0$. 

These two empirical observations --  that random $k$-SAT exhibits a sharp threshold  and that instances near the threshold are computationally hard -- set off an explosion of work on random $k$-SAT  and related models (random graph coloring, $k$-NAE-SAT, $k$-XOR-SAT, etc.) in many different directions.  

We can model random $k$-SAT in the setting of this paper.  For a given $n$ and $p$, let $F_k(n,p)$ be a random $k$-SAT formula generated by including each of the $\binom{2n}{k}$ possible $k$-clauses independently with probability $p$.  The property of being satisfiable is a non-trivial property with a monotone complement (and thus has a threshold function).

One long-standing conjecture is that the random $k$-SAT model exhibits a sharp threshold.
\begin{conj}[Satisfiability Conjecture]
\label{conjSAT}
For each $k \ge 2$, there exists $c_k >0$ so that for every $\eps>0$, the following hold:
\begin{itemize}
\item If $p \le (1-\eps) c_k n^{-1/(k-1)}$, then whp $F_k(n,p)$ is satisfiable.
\item If $p \ge (1+\eps) c_k n^{-1/(k-1)}$, then whp $F_k(n,p)$ is unsatisfiable.
\end{itemize}
\end{conj}
Conjecture~\ref{conjSAT} was proved early on for the special case of random $2$-SAT~\cite{chvatal1992mick}; and in fact the precise scaling window was determined~\cite{bollobas2001scaling}. The reason for this is that determining satisfiability of a $2$-SAT formula can be reduced to determining  the strongly connected components of the implication graph on literals; this gives a linear-time algorithm for $2$-SAT in the worst case and gives a strong analogy between the random $2$-SAT threshold and the emergence of a giant component in a random graph. The case $k\ge 3$ is fundamentally different.

Random $k$-SAT is an example of a \textit{random constraint satisfaction problem} (random CSP); there is an initial  set of possible solutions (here all possible assignments to $n$ boolean variables) and a set of random constraints is selected, each of which rules out some solutions (here each clause rules out a $2^{1-k}$ fraction of all possible solutions).  Random graph coloring is another random CSP: the $q^n$  colorings of $n$ vertices are the possible solutions; each edge rules out $\frac{1}{q}$-fraction of the solutions.

How can one locate the satisfiability threshold in a random CSP?  As with many of the examples above, a starting point is to try the first- and second-moment methods on the random variable $Z$ that counts the number of solutions.  If $p$ is large enough that $\E Z \to 0$, we know $p_c \le p$; conversely if $p$ is such that $\E Z \to \infty$ and, say, $\E (Z^2) \le C (\E Z )^2$, then (via the Paley--Zygmund Inequality and Friedgut's Theorem, see Sections~\ref{secFried} and~\ref{secCavity}) we know $p_c \ge p$. 

This approach, combined with smart ideas of conditioning, replacing $Z$ with a related random variable, and solving difficult optimization problems to bound $\E (Z^2)$, has been able to pin down the satisfiability threshold in many random CSP's to within a small constant factor (e.g.~\cite{kirousis1998approximating,achlioptas2002asymptotic,achlioptas2003threshold,achlioptas2005rigorous}). These probabilistic methods face a barrier, however, and cannot address Conjecture~\ref{conjSAT} or determine $p_c$.  Instead, new tools and ideas from mathematics and statistical physics have been brought to bear on the problem.

\subsection{Friedgut's Theorem}
\label{secFried}

Major progress towards Conjecture~\ref{conjSAT} was made by Friedgut in~\cite{friedgut1999sharp} who proved an $n$-dependent sharp threshold.
\begin{theorem}[Friedgut]
\label{thmFksat}
For each $k \ge 3$ there is a function $c_k(n)$ bounded above and below by constants so that for every $\eps>0$ the following hold:
\begin{itemize}
\item If $p \le (1-\eps) c_k(n) \cdot n^{-1/(k-1)}$, then whp $F_k(n,p)$ is satisfiable.
\item If $p \ge (1+\eps) c_k(n)  \cdot n^{-1/(k-1)}$, then whp $F_k(n,p)$ is unsatisfiable.
\end{itemize}
\end{theorem}
This theorem shows that the scaling window of random $k$-SAT is $o(p_c)$ but it leaves open the possibility that the critical density $c_k$ fluctuates as $n \to \infty$.  

Though motivated by the Satisfiability Conjecture, Friedgut proved much more in~\cite{friedgut1999sharp}: he gave a very general characterization of which monotone properties of random graphs can have a coarse threshold: those that are well approximated by the property of containing some bounded-size subgraph from a bounded-size list of subgraphs.  All other properties have a sharp threshold (in the non-uniform sense of Theorem~\ref{thmFksat}).  

Returning to the example from Section~\ref{sec3examples}, Achlioptas and Friedgut~\cite{achlioptas1999sharp} use~\cite{friedgut1999sharp} to show that for $q \ge 3$ the property of $G(n,p)$ being $q$-colorable exhibits a sharp threshold, in the sense of Theorem~\ref{thmFksat}.  It remains open to show that there exists constants $d_q$ so that a sharp threshold for $q$-colorability occurs at $d_q/n$.

\subsection{The cavity method and the structure of  solution spaces}
\label{secCavity}

Statistical physicists soon turned their attention to random $k$-SAT and related models using tools and intuition from the study of spin glasses~\cite{mezard1987spin,monasson1997statistical,mezard2002analytic,mertens2006threshold,krzakala2007gibbs,montanari2008clusters}.  The replica and cavity methods they developed are  powerful analytic tools for analyzing disordered (random) systems, based on assumptions that to a large extent have not been proved mathematically.  Both are methods for computing $\E \log Z$, the expectation of a log partition function, the normalizing constant of a Gibbs measure\footnote{Here the setting is always at positive temperature so $Z>0$ and taking a logarithm is justified; models with hard constraints like random $k$-SAT or $q$-coloring are analyzed by taking a limit as the temperature tends to $0$.}.  The replica method is based on the identity $\log Z = \lim_{n \to 0} \frac{Z^n -1}{n}$, and rather magically using computations for large integer $n$ to deduce the value of the limit.  See~\cite[Chapter 1.13]{talagrand2010mean} for a mathematical introduction. The cavity method is based on a coupling of random graphs or hypergraphs on $n$ vertices and on $n+1$ vertices, and understanding the expected change in the partition function, $\log Z_{n+1} - \log Z_n$, as a function of the Gibbs measure   and its correlation properties.  See~\cite[Chapter 19]{mezard2009information} for an introduction to the method. Under  strong assumptions, these two methods are capable of making very detailed predictions for thresholds and phase transitions in a broad class of random structures, \textit{Gibbs measures on random graphs}.  

Applying these methods to random CSP's involves shifting perspective and instead of asking whether or not there is a solution, asking `what does the uniform distribution on solutions (if they exist) look like?'  Then the main properties of interest are \textit{correlations} and \textit{overlaps}.  Correlations are measured, for instance, by the covariance in the indicator random variables that two variables take the value True (or two vertices take the color Red).  Overlap is the random variable measuring the fraction of coordinates on which two independent samples from the uniform distribution on solutions agree.

\begin{figure}
\begin{centering}
\includegraphics[width=.78\linewidth]{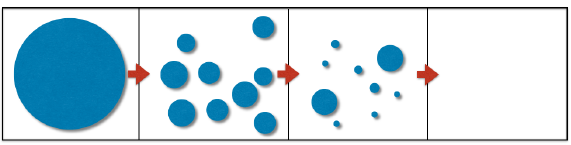}
\caption{Cartoon of the evolution of the solution space for random $k$-SAT or random graph coloring as the edge density increases.}
\label{figCavity}
\end{centering}
\end{figure}

Figure~\ref{figCavity} depicts a cartoon of how the space of solutions  in, say, random $k$-SAT or random graph $q$-coloring, changes (or is predicted to change) as the density of constraints increases.  Three distinct thresholds are pictured.  At low densities, solutions lie in one large component, connected under single variable changes.  Then after the \textit{shattering threshold} solutions break apart into exponentially many clusters of roughly equal exponential size, each separated by linear Hamming distance ($\Theta(n)$ variables must be changed to move between clusters). Next, at the \textit{condensation threshold} a constant number of large clusters contain almost all solutions, while there are exponentially many smaller clusters. Finally, after the \textit{satisfiability threshold}, no solutions remain.

The condensation threshold is critical for correlations and overlap.  In the random graph coloring model, below the condensation threshold average correlations between vertices vanish as $n\to \infty$, and the overlap concentrates on $1/q$ (as it would if the graph were empty and we were sampling uniformly from $[q]^n$). Above the threshold, average correlations are bounded away from $0$ and overlap concentrates on two points, $1/q$ and some $\eta > 1/q$, with the second corresponding to the case in which both samples come from the same dominant cluster.  In the language of the cavity method, the model is \textit{replica symmetric} below the condensation threshold and exhibits \textit{replica symmetry breaking} above the threshold.

In terms of identifying the satisfiability threshold, condensation is important because it presents a fundamental obstacle to applying the second-moment method to the random variable $Z$ counting the number of solutions.  We briefly describe the connection at a high level.   Let $d$ denote the average vertex or variable degree of a random CSP.  The exponential order of the first moment, $f_1(d) := \lim_{n \to \infty} \frac{1}{n} \log \E Z$, is a decreasing function of $d$ and the first-moment threshold $d_{\mathrm{first}}$ is the point at which it hits $0$.  For $d>d_{\mathrm{first}}$ whp there are no solutions by Markov's Inequality; for $d< d_{\mathrm{first}}$ the expected number of solutions is exponentially large. To show that there is at least one solution with probability bounded away from $0$, we would like to apply the second-moment method via the Paley--Zygmund Inequality: $\Pr(Z >0 ) \ge   \frac{ ( \E Z)^2}{\E (Z^2)  }$.  Writing $\E(Z^2) = \sum_{\sigma, \tau} \Pr[ \sigma, \tau \text{ both solutions}]$ one sees that on an exponential scale, the second moment is at most the square of the first moment if and only if the dominant contribution to $\E(Z^2)$ comes from pairs $\sigma, \tau$ of uncorrelated vectors.   Past the condensation threshold, pairs of solutions with non-trivial correlation are significant and so the second moment should be expected to be exponentially larger than the square of the first moment. See the discussion in~\cite{achlioptas2006random,coja2012catching,coja2012condensation,coja2014asymptotic,bapst2016condensation} for much more detail.  To prove a random CSP is whp satisfiable beyond the condensation threshold requires a significantly different approach.    

\subsection{Proof of the satisfiability conjecture for large $k$}

How can one locate the satisfiability threshold in light of condensation blocking the second-moment method?  One very successful solution is to use the cavity method predictions to select a different random variable on which to perform the second-moment method.  In particular, the predicted `1-RSB' behavior of both random $k$-SAT and random graph coloring suggests that while the uniform distribution over solutions exhibits replica symmetry breaking above the condensation threshold, the uniform distribution on \textit{clusters} does not. 

The innovation of Coja-Oghlan and Panagiotou in~\cite{coja2012catching} was to design (for the closely related random $k$-NAE SAT model\footnote{`not-all-equal' SAT: an assignment is satisfying if every clause has at least one true and at least one false literal.}) a random variable $Z_\beta$ counting the number of solutions with $\beta n$ `blocked' variables that is a stand-in for counting the number of clusters of a given (exponential) size.  The same authors apply this idea to random $k$-SAT in~\cite{coja2014asymptotic}, additionally conditioning on typical vertex degrees.  This kind of random variable was also used in~\cite{ding2016satisfiability,ding2016maximum} to determine thresholds and limiting values of optimization problems in random regular models (in which each vertex or variable has the same edge or constraint degree).

In~\cite{ding2022proof}, Ding, Sly, and Sun used this cavity-method inspired second-moment argument along with conditioning on the typical empirical distribution of arbitrary depth neighborhoods of variables to resolve Conjecture~\ref{conjSAT} for $k$ large enough.  
\begin{theorem}[Ding, Sly, Sun~\ref{conjSAT}]
There exists $k_0 >0$ so that for each $k \ge k_0$ there exists $c_k$ so that the following hold. 
\begin{itemize}
\item If $m \le (1-\eps) c_k n$, then whp $F_k(n,m)$ is satisfiable.
\item If $m \ge (1+\eps) c_k n$, then whp $F_k(n,m)$ is unsatisfiable.
\end{itemize}
\end{theorem}
The proof is a mathematical tour-de-force, employing several sophisticated tools (for both the upper and lower bounds) and performing very challenging probabilistic and combinatorial calculations.  Notably, to prove Conjecture~\ref{conjSAT} in this case they establish the exact threshold as predicted by the 1-RSB cavity method. 

One could hope to apply the same strategy to random graph coloring, but this seems hopelessly complex: having $q >2$ possibilities for each variable makes everything more complicated, and the proof in~\cite{ding2022proof} is already over 300 pages.

\subsection{Algorithmic thresholds}

What about the other question asked by Mitchell, Selman, and Levesque --  is it true that random $k$-SAT instances near the satisfiability threshold are computationally intractable? 

There was some hope among statistical physicists that \textit{survey propagation}, a message-passing algorithm based on the 1-RSB predictions for random $k$-SAT and other models, might provide an efficient search algorithm for instances near the satisfiability threshold~\cite{braunstein2005survey}. However (at least for large $k$) this was disproved~\cite{hetterich2016analysing,bresler2022algorithmic}.

Instead, evidence of computational hardness based on the structure of the solution space has emerged.  Achlioptas and Coja-Oghlan proved the shattering of the solution space in~\cite{achlioptas2008algorithmic} and observed that the threshold for shattering approximately coincides with the density above which no efficient search algorithms are known (though see below for a caveat in linking the two).     While it did not establish a direct link between solution space structure and algorithms, this paper was innovative both in its techniques (using the \textit{planted model}, described more below) and in making a conceptual link between the two.

A more direct link with algorithms is the Overlap Gap Property (OGP) pioneered by Gamarnik and Sudan~\cite{gamarnik2017limits,gamarnik2021overlap}.  The basic OGP states that there is an interval $(a,b)$ so that whp over an instance of a random computational problem, no pair of solutions have their normalized overlap in $(a,b)$.  This can be proved using a first-moment method.  Amazingly, this simple property then implies that entire classes of algorithms (local algorithms, low-degree algorithms) cannot find a solution whp.  For random $k$-SAT the OGP perspective has been applied to determine approximately the low-degree algorithmic threshold of the problem~\cite{bresler2022algorithmic}.

\subsection{Scaling windows}

With Friedgut's theorem establishing a sharp threshold in the random $k$-SAT model, one could hope to say something quantitatively stronger about scaling windows in this and related models.    Friedgut's theorem (or Bourgain's more abstract theorem) gives bounds on the $k$-SAT scaling window of the form $T_\eps = O( p_c / \log \log n)$, just barely enough for a sharp threshold.  Surprisingly, Abbe and Montanari~\cite{abbe2014concentration} show that even a very mild improvement to this bound would resolve Conjecture~\ref{conjSAT}  
\begin{prop}[\cite{abbe2014concentration}]
If for any fixed $\del>0$, the scaling window for random $k$-SAT satisfies $T_\eps = O( p_c /\log^{1+\del} n)$, then Conjecture~\ref{conjSAT} holds.
\end{prop}

The statistical physics cavity method does not give much guidance on what to expect for scaling windows.  In 2002, Wilson~\cite{wilson2002critical} used basic probabilistic arguments to prove generic polynomial lower bounds on the scaling window for random $k$-SAT and other models, thus disproving a number of conjectures from statistical physics.  Good upper bounds for random $k$-SAT or random graph $q$-coloring are completely lacking currently.

In a different random CSP, however, scaling windows have recently been pinned down precisely.   The \textit{binary perceptron}~\cite{cover1965geometrical} arose in the 1960's as a toy model of a neural network then attracted attention in statistical physics via the work of Gardner and others~\cite{G87,GD88,KM89}.  A symmetric variant was recently introduced in~\cite{aubin2019storage} and has some remarkably nice properties from a mathematical point of view.  For one, the plain first- and second-moment methods (along with some concentration arguments) suffice to pin down precisely the satisfiability threshold in the model~\cite{aubin2019storage,perkins2021frozen,abbe2022proof}.  Closely related is the fact that the solution space looks very different than the cartoon in Figure~\ref{figCavity}: for all constraint densities below the satisfiability threshold, whp over the instance, almost all solutions are isolated, at linear Hamming distance from the nearest other solution.  The corresponding cartoon would be a sprinkling of points.  

The nature of the solution space in this model raises a lot of questions: if almost all solutions are isolated should the search problem be hard at all densities?  In fact there are efficient search algorithms at low densities~\cite{kim1998covering,abbe2022binary} even when almost all solutions are isolated; these algorithms in fact find solutions lying in  very rare clusters (with maximum possible diameter)~\cite{abbe2022binary}, as predicted by physicists working on machine learning problems~\cite{baldassi2015subdominant}.  These algorithms are  `multiscale majority algorithms' in which entries of a potential solution vector $X$ are set round by round to optimize a greedy objective. The symmetric binary perceptron is a random model of the combinatorial discrepancy problem, and so constructive discrepancy minimization algorithms also provide provably efficient search algorithms at low constraint densities~\cite{lovett2015constructive,rothvoss2017constructive,eldan2018efficient,bansal2020line}.  These algorithms can also be adapted for the original (non-symmetric) perceptron model~\cite{li2024discrepancy}.  The success of these algorithms indicates that the cartoon in Figure~\ref{figCavity} is not really relevant to algorithms: the cartoon depicts properties of typical solutions, while efficient algorithms may indeed find rare solutions.  The structural properties like the OGP that apply to \textit{all} solutions (or tuples of solutions) are more algorithmically relevant. 

Finally, in a recent breakthrough Altschuler proved something remarkable about the scaling window: the scaling window, measured in number of constraints, is of width $O(\log n)$ (while the critical number of constraints is $O(n)$)~\cite{altschuler2022critical}.  Contrast this to what is known via abstract results like Friedgut's theorem which would give a bound  of $O(n/\log \log n)$; and to what is known for models like random $k$-SAT and random graph coloring which is only what the abstract results give.  Even more recently Sah and Sawhney determined the scaling window completely, showing it is of width $O(1)$ and giving the limiting probability of satisfiability inside the window~\cite{sah2023distribution}.

\section{Statistical inference and the stochastic block model}
\label{secStat}

A fundamental statistical question can be phrased as `Under what circumstances can a signal be recovered from a noisy observation of that signal?'  Or - `Can we distinguish a signal from noise?'    

As the study of statistics has evolved in the age of fast computers and massive data sets, these same questions remain fundamental but the data sets of interest are now very high dimensional and the question of efficient computation becomes paramount.   Just as we want to know whether `typical' instances of $3$-SAT are algorithmically tractable we also want to know if `typical' statistical inference problems are algorithmically tractable.  

A very useful framework for the rigorous study of these questions is the \textit{teacher--student} framework (see e.g.,~\cite{zdeborova2016statistical}).  In this framework, a teacher describes a generative probabilistic model of data to a student. The model takes as an input a `ground truth' and adds to this noise of some form.  The teacher chooses a ground truth from some known prior distribution, generates the data from the model, and presents the data to the student.  The student's task is to recover the ground truth from the data, using knowledge of the generative model and prior distribution.

Perhaps the most studied teacher--student model in the \textit{stochastic block model}, a toy model for the statistical and machine learning task of clustering: partitioning a set of data points into subsets with similar characteristics.   The stochastic block model deals with a simple specialization of clustering: community detection in which the the data is a graph and the task is to partition the vertex set into subsets of vertices with similar connectivity structure; this could mean finding a partition in which most edges lie within parts, or finding a partition in which most edges cross the parts. 

Formally, a symmetric version of the stochastic block model is defined as follows.  Fix an integer $q \ge 2$, and $p_{\text{in}}, p_{\text{out}} \in [0,1]$, $p\ne q$, and let $n$ denote the number of vertices of the graph.
\begin{itemize}
\item Choose a partition  $\sigma \in [q]^n$  of $[n]$ into $q$ parts uniformly at random.
\item For each pair $i, j \in [n]$, $i \ne j$, include the edge $(i,j)$ with probability $p_{\text{in}}$ if $\sigma(i) = \sigma(j)$ and with probability $p_{\text{out}}$ if $\sigma(i) \ne \sigma(j)$, all edges independently of the others. 
\item Call the resulting graph $\hat G$. 
\end{itemize}

If $p_{\text{in}} >p_{\text{out}}$, then, on average, more edges will be drawn between vertices with the same label (we say the model is \textit{assortative}); if $p_{\text{in}} >p_{\text{out}}$, the opposite is true (the model is \textit{disassortative}).   If $p_{\text{in}} =p_{\text{out}}$ the model is simply $G(n,p)$.

The inference task is to \textit{recover} the partition $\sigma$ given the graph $\hat G$.  Recovery can be defined in different ways: exact recovery (up to a permutation of the $q$ labels); recovery of almost all the labels up to permutation; or recovery of a partition $\sigma'$ that (after a permutation) agrees with $\sigma$ on $\frac{1}{q} + \eps$ fraction of vertices, for some small $\eps >0$; that is, just a tiny bit better than random guessing.   We focus here on the last notion, called \textit{weak recovery}.

Though we have left the setting of monotone properties of the hypercube under $\mu_p$, we have retained some of the essential features: this model has independent edges and an important monotonicity property.  If we fix the ratio $\frac{p_{\text{in}}}{p_{\text{out}}} =: \alpha$ and set $p_{\text{in}}= p \alpha $ and $p_{\text{out}}= p$  then the weak recovery problem only becomes \textit{easier} as $p$ increases: we have more data (observed edges) from which to deduce the signal (the underlying partition).  That is, if $p > p'$ we can reduce the model with parameter $p$ to that with parameter $p'$ by independently deleting each edge with probability $p'/p$.

\subsection{Sharp thresholds for inference}

Though the stochastic block model was  defined independently in different fields in the 1980's~\cite{holland1983stochastic,bui1987graph,dyer1989solution}, mathematical interest in the model exploded in 2011 when Decelle,  Krzakala, Moore, and Zdeborov{\'a}~\cite{decelle2011asymptotic} used the cavity method to make a series of beautiful predictions about sharp thresholds in the model.

We focus on their conjectures for thresholds  for weak recovery when $p_{\text{in}}, p_{\text{out}} = O(1/n)$ and so the random graph has constant average degree.  

Fix a number of parts $q \ge 2$ and a ratio $\alpha>0$ of $p_{\text{in}}/p_{\text{out}}$. 
There is an  \textit{information theoretic threshold} at $d_{\text{inf}}$ if 
\begin{enumerate}
\item when $d < d_{\text{inf}}$, $p_{\text{in}}= \alpha \frac{d}{n}$, $p_{\text{out}}= \frac{d}{n}$, there is no algorithm (efficient or not) that whp finds a partition $\sigma'$ that agrees with $\sigma$  on $\frac{1}{q} + \eps$ fraction of vertices (after a permutation of the parts) for any fixed $\eps>0$.
\item when $d > d_{\text{inf}}$, $p_{\text{in}}= \alpha \frac{d}{n}$, $p_{\text{out}}= \frac{d}{n}$, there is $\eps>0$ and an  algorithm (perhaps inefficient) that whp finds a partition $\sigma'$ that agrees with $\sigma$  on $\frac{1}{q} + \eps$ fraction of vertices (after a permutation of the parts).
\end{enumerate}

There is an  \textit{algorithmic threshold} at $d_{\text{alg}}$ if 
\begin{enumerate}
\item when $d < d_{\text{alg}}$, $p_{\text{in}}= \alpha \frac{d}{n}$, $p_{\text{out}}= \frac{d}{n}$, there is no polynomial-time algorithm that whp finds a partition $\sigma'$ that agrees with $\sigma$  on $\frac{1}{q} + \eps$ fraction of vertices (after a permutation of the parts) for any fixed $\eps>0$.
\item when $d > d_{\text{alg}}$, $p_{\text{in}}= \alpha \frac{d}{n}$, $p_{\text{out}}= \frac{d}{n}$, there is $\eps>0$ and a polynomial-time algorithm  that whp finds a partition $\sigma'$ that agrees with $\sigma$  on $\frac{1}{q} + \eps$ fraction of vertices (after a permutation of the parts).
\end{enumerate}

For the case $q=2,3$ they conjectured that for any $\alpha$,  $d_{\text{inf}}=d_{\text{alg}}$, but for $q \ge 5$ they conjectured a gap: $d_{\text{inf}}<d_{\text{alg}}$ (now known as a statistical--computational gap, and the subject of great recent interest in statistics and computer science).  In general the conjectured value $d_{\text{alg}}$ is explicit (it is the `Kestum--Stigum threshold'~\cite{kesten1966limit}), while the conjectured value of $d_{\text{inf}}$ is given implicitly as the solution of a complicated optimization problem over probability measures.    The $q=2$ conjectures were proved in celebrated works of Mossel--Neeman--Sly~\cite{mossel2015reconstruction,mossel2018proof} and Massouli{\'e}~\cite{massoulie2014community}.  

For $q \ge 3$, partial progress has been made.  
Abbe and Sandon~\cite{abbe2018proof} proved that for any $q$ and any $\alpha$, weak recovery is possible above the Kestum--Stigum threshold, i.e. the positive part of the conjecture on the algorithmic threshold.
 For $q \ge 3$ and $\alpha <1$ (the disassortative case), the conjecture on the value of the information--theoretic threshold (as the solution of the complicated optimization problem) was proved in~\cite{coja2018information} by implementing a rigorous version of the cavity method.  Recently, it was proved in~\cite{mossel2023exact} that for $q=3$ (and for $q=4$ and $\alpha$ sufficiently close to $1$), $d_{\text{inf}}=d_{\text{alg}}$.

\subsection{The planted model}
\label{secPlanted}

To wrap things up, we discuss a concept with connections to all of the developments recounted in Sections~\ref{secLatin},~\ref{secCSP}, and~\ref{secStat}.  

The stochastic block model is an example of a \textit{planted model}.  A solution (in this case a partition into $q$ parts) is chosen, then a random instance is drawn consistent with this \textit{planted solution}.    The inference task above is to recover information about the planted solution given the random instance.  An even easier computation task is to distinguish the planted model from the purely random model (in this case the \ER random graph $G(n,p)$).

It is straightforward to devise planted models for random $k$-SAT and random graph coloring: pick a solution $\sigma$ uniformly at random and choose independent constraints or edges among those that are satisfied by $\sigma$.

Achlioptas and Coja-Oghlan~\cite{achlioptas2008algorithmic} show that if the number of solutions in the random model is sufficiently concentrated, then high probability results about the planted solution in the planted model can be transferred to high probability statements about uniformly random solutions in the random model.  This is the key to establishing results on the structure of the solution space (and is used in the symmetric perceptron results~\cite{perkins2021frozen,abbe2022proof} as well).

Returning again to random graph coloring, the planted model is exactly the extreme case ($\alpha=0$) of the disassortative stochastic block model.   In~\cite{bapst2016condensation,coja2018information,coja2018charting}  the condensation threshold for random graph $q$-coloring is determined precisely, in terms of a solution to a variational problem arising from the cavity method.  This bound is the best known lower bound on the $q$-colorability threshold for $q\ge 4$ (the best bound for $q=3$ is in~\cite{achlioptas2002almost}).  In~\cite{coja2018information}, a precise connection between condensation thresholds and information theoretic thresholds is made: the two thresholds coincide for a large class of models exhibiting some symmetry and convexity; these include random graph coloring, the anti-ferromagnetic Potts model, $k$-NAE-SAT. They do not include random $k$-SAT (asymmetry) or the assortative stochastic block model (lacking the needed convexity).  Overcoming these technical obstacles and being able to identify condensation thresholds in asymmetric models is a major challenge and could eventually lead to a generic approach to Conjecture~\ref{conjSAT}.  See related discussion in~\cite{yu2023ising,kireeva2023breakdown}.

To conclude, somewhat surprisingly the planted model is closely related to the developments recounted in Section~\ref{secLatin} as well.
Mossel, Niles-Weed, Sun and Zadik~\cite{mossel2022second} recently showed that a well-chosen planted model and a second-moment argument can be used to prove the `spread lemma', the key technical ingredient in the proof of improved sunflower bounds~\cite{alweiss2021improved} and the fractional Kahn--Kalai conjecture~\cite{frankston2021thresholds}. 

\section{Conclusions and questions}
\label{secConclude}

As we saw in Section~\ref{secLatin}, the fractional Kahn--Kalai theorem has sufficed for all known applications, but one can ask if there are applications that require Theorem~\ref{thmKK}.  

\begin{question}
Are there threshold applications that need the full power of Theorem~\ref{thmKK} rather than the fractional version of Theorem~\ref{thmFrac}? Or is there a constant $K$ so that $q_f(\cF) \le K q(\cF)$? 
\end{question}
The second is conjectured in the affirmative by Talagrand in~\cite{talagrand2010many}.

As we have seen, one effective way of determining sharp thresholds and scaling windows is to prove a hitting time result, relating a complex property to a simple property.  

\begin{question}
Is there a general way to find first-order asymptotics of $p_c$ for the property of $G(n,p)$ containing a subgraph $H$ when there is not a natural hitting-time conjecture?
\end{question}

A nice conjecture for one such property is by Kahn, Narayanan, and Park on the threshold for the existence of the square of a Hamilton cycles in $G(n,p)$; that is, the existence of a cyclic ordering of the $n$ vertices so that there is an edge between every pair of neighbors and second neighbors.
\begin{conj}[Kahn, Narayanan, and Park~\cite{kahn2021threshold}]
For the property $\cF$ of $G(n,p)$ containing the square of a Hamilton cycle, 
\[ p_c(\cF) = (1+o(1)) \sqrt{\frac{e}{n}} \,. \]
\end{conj}
In~\cite{kahn2021threshold} they prove $p_c = \Theta (n^{-1/2})$ using a finer understanding of the spreadness of the uniform distribution on squares of Hamilton cycles; a direct application of Theorem~\ref{thmSpread} would have lost a $\log n$ factor.  See also~\cite{spiro2023smoother,espuny2023} for  generalizations.

\begin{question}
\label{qScalingTools}
Are there any general-purpose tools for bounding the width of the scaling window of a monotone property with a sharp threshold, beyond the bounds given by Friedgut's theorem?  Can Conjecture~\ref{conjSAT} be proved in a generic way, without the need to precisely identify the threshold?
\end{question}

\subsection*{Random geometric graphs}

We conclude the survey by stating one of the author's favorite open problems on thresholds in random structures.

In this survey we have discussed the \ER random graph in depth and another model of a random graph in the stochastic block model; both models have the important property that edges are independent.  There are many other important random graph models that do not have this property.  Here we discuss one, the \textit{random geometric graph}~\cite{gilbert1961random,penrose2003random}, in which edges between randomly placed points are determined geometrically.

In particular, let $\mathbb S^{d-1}$ be the unit sphere in $d$-dimensions\footnote{We alternatively could  choose the underlying space to be the $d$-dimensional unit torus, $\mathbb T^d =  \R^d/\Z^d$, the results and questions below would still apply.}.  The random graph $G_{d}(n,p) = (V,E)$ is formed as follows:
\begin{itemize}
\item Let $V = [n]$
\item Select $n$ points $x_1, \dots, x_n$ independently and uniformly from $\mathbb S^{d-1}$.  
\item Let $E = \{ (i,j) : x_i \cdot x_j \ge \tau_p \}$ where $\tau_p$ is chosen so that $\Pr[x_i \cdot x_j \ge \tau_p] =p$.
\end{itemize}
In particular, like the \ER random graph $G(n,p)$,  $G_d(n,p)$ is a random graph on $n$ vertices with edge probability $p$; but now the edges are not independent. 

Thresholds in a random geometric graphs have been studied extensively; see the textbook of Penrose~\cite{penrose2003random} for results on connectivity, subgraph containment, existence of a giant component.  More recently, following the influential paper of Devroye, Gy{\"o}rgy, Lugosi, and Udina~\cite{devroye2011high}, thresholds in $G_d(n,p)$ when $d= d(n) \to \infty$ have been studied, along with a different kind of threshold: for a given $p(n)$ for which values of $d(n)$ are the two random graphs $G(n,p)$ and $G_d(n,p)$ statistically distinguishable?~\cite{bubeck2016testing,liu2022testing}.

Just as in $G(n,p)$, a general study of monotone properties is possible in $G_d(n,p)$. In particular, let $\cF$ be a  non-trivial monotone  property (where, as always, monotone pertains to adding edges).  Then again the probability that $\cF$ holds in $G_d(n,p)$ is a strictly increasing function of $p$, and so $p_c(\cF)$ can be uniquely defined, and all of the same questions from 
Section~\ref{secProperties} can be asked in this setting too.

Perhaps surprisingly then, the analogue of Theorem~\ref{thmBT}, the Bollob{\'a}s-Thomason theorem, that lays the foundation for the general study of thresholds of monotone properties in random graphs, is not known in general for random geometric graphs.  In the special case of $d=1$ (points on a circle), McColm proved  that every monotone property has a threshold~\cite{mccolm2004threshold}.  Moreover, Goel, Rai, and Krishnamachari~\cite{goel2005monotone} proved that the scaling window of every monotone property is bounded by a function (depending on $d$) that vanishes as $n \to \infty$, the analogue of the result of Friedgut and Kalai~\cite{friedgut1996every} for $G(n,p)$.  The general statement, however, remains open.
\begin{conj}
\label{conjRGG}
Every non-trivial monotone property has a threshold in $G_d(n,p)$.  
\end{conj}

\section*{Acknowledgements}

The author thanks Vishesh Jain, Tom Kelly, and Richard Montgomery for very helpful comments and discussions while preparing this survey.

\bibliography{thresholds.bib}
\bibliographystyle{abbrv}

\end{document}